\documentclass[11pt]{article}
\usepackage{amsfonts}
\usepackage{amssymb}
\usepackage{amsmath}

\def\sqr#1#2{{\vcenter{\vbox{\hrule height.#2pt
              \hbox{\vrule width.#2pt height#1pt \kern#1pt \vrule
width.#2pt}
              \hrule height.#2pt}}}}
\def\signed #1{{\unskip\nobreak\hfil\penalty50
              \hskip2em\hbox{}\nobreak\hfil#1
              \parfillskip=0pt \finalhyphendemerits=0 \par}}
\def\endpf{\signed {$\sqr69$}}
\def\dbN{{\mathbb{N}}}

\def\dbR{{\mathbb{R}}}

\def\e{\varepsilon}

\def\3n{\negthinspace \negthinspace \negthinspace }
\def\2n{\negthinspace \negthinspace }
\def\1n{\negthinspace }
\def\ns{\noalign{\smallskip} }
\def\ns{\noalign{\medskip} }
\def\ds{\displaystyle}
\def\O{\Omega}
\def\no{\noindent}

\def\ms{\medskip}

\def\q{\quad}
\def\qq{\qquad}
\def\hb{\hbox}

\def\esssup{\mathop{\rm esssup}}

\def\({\Big (}
\def\){\Big )}
\def\[{\Big[}
\def\]{\Big]}

\def\be{\begin{equation}}
\def\bel{\begin{equation}\label}
\def\ee{\end{equation}}
\def\bea{\begin{eqnarray}}
\def\eea{\end{eqnarray}}
\def\bt{\begin{theorem}}
\def\et{\end{theorem}}
\def\bc{\begin{corollary}}
\def\ec{\end{corollary}}
\def\bl{\begin{lemma}}
\def\el{\end{lemma}}
\def\bp{\begin{proposition}}
\def\ep{\end{proposition}}
\def\br{\begin{remark}}
\def\er{\end{remark}}
\def\ba{\begin{array}}
\def\ea{\end{array}}
\def\bd{\begin{definition}}
\def\ed{\end{definition}}

\newtheorem{lemma}{Lemma}[section]
\newtheorem{remark}{Remark}[section]

\newtheorem{theorem}{Theorem}[section]
\newtheorem{corollary}{Corollary}[section]

\newtheorem{definition}{Definition}[section]
\newtheorem{proposition}{Proposition}[section]

\makeatletter
   
   \@addtoreset{equation}{section}
\makeatother

\begin{document}

\title{\bf Weak Maximum Principle for Strongly Coupled
Elliptic Differential Systems}
\author{Xu Liu\thanks{School of Mathematics and Statistics, Northeast Normal
University, Changchun 130024, China. This work is supported by the
NSF of China under grant 10901032. E-mail address:
liuxu@amss.ac.cn.\ms}\ \ \ and\ \  \ Xu Zhang\thanks{Key Laboratory
of Systems and Control, Academy of Mathematics and Systems Science,
Chinese Academy of Sciences, Beijing 100190, China; Yangtze Center
of Mathematics, Sichuan University, Chengdu 610064, China. This work
is supported by the NSF of China under grant 10831007 and by the
National Basic Research Program of China (973 Program) under grant
2011CB808002. E-mail address: xuzhang@amss.ac.cn.}}

\date{}

\maketitle
\begin{abstract}

A classical counterexample due to E. De Giorgi, shows that the weak
maximum principle does not remain true for general linear elliptic
differential systems. After that, there are some efforts to
establish the weak maximum principle for special elliptic
differential systems, but the existing works are addressing only the
cases of weakly coupled systems, or almost-diagonal systems, or even
some systems coupling in various lower order terms.  In this paper,
by contrast, we present maximum modulus estimates for weak solutions
to two classes of coupled linear elliptic differential systems with
different principal parts, under considerably mild and physically
reasonable assumptions. The systems under consideration are strongly
coupled in the second order terms and other lower order terms,
without restrictions on the size of ratios of the different
principal part coefficients, or on the number of equations and space
variables.
\end{abstract}



\section{Introduction}

Let $m,n\in \dbN\backslash\{0\}$, and $\Omega\subset \dbR^m$ be a
bounded domain with an $C^1$ boundary $\Gamma$ and having the cone
property. We consider the following nonhomogeneous, isotropic
elliptic differential system of second order:
\begin{eqnarray}\label{1}
\left\{
\begin{array}{lr}
-\mbox{div}(a^{11}\nabla y^1)-\mbox{div}(a^{12}\nabla
y^2)-\cdots-\mbox{div}(a^{1n}\nabla
y^n)+\displaystyle\sum_{i=1}^nC^{1i}\cdot\nabla y^i+D^1\cdot y=f^1
&\mbox{ in }\Omega,\\\ns -\mbox{div}(a^{21}\nabla
y^1)-\mbox{div}(a^{22}\nabla y^2)-\cdots-\mbox{div}(a^{2n}\nabla
y^n)+\displaystyle\sum_{i=1}^nC^{2i}\cdot\nabla y^i+D^2\cdot y=f^2
&\mbox{ in }\Omega,\\\ns
\quad\quad\quad\quad\quad\quad\quad\quad\quad
\quad\quad\quad\quad\quad\quad\quad\vdots\\\ns
-\mbox{div}(a^{n1}\nabla y^1)-\mbox{div}(a^{n2}\nabla
y^2)-\cdots-\mbox{div}(a^{nn}\nabla
y^n)+\displaystyle\sum_{i=1}^nC^{ni}\cdot\nabla y^i+D^n\cdot y=f^n
&\mbox{ in }\Omega,\\\ns y^1=g^1,\quad y^2=g^2,\quad \cdots,\quad
y^n=g^n  &\mbox{ on }\Gamma,
\end{array}
\right.
\end{eqnarray}
and the following general nonhomogeneous elliptic differential
system of second order:
\begin{eqnarray}\label{61}\left\{
\begin{array}{lr}
-\displaystyle\sum_{p, q=1}^{m}\left[(a_{p
q}^{11}y^1_{x_p})_{x_q}+(a_{p q}^{12}y^2_{x_p})_{x_q}+\cdots+(a_{p
q}^{1n}y^n_{x_p})_{x_q}\right]+\displaystyle\sum_{i=1}^nC^{1i}\cdot\nabla
y^i+D^1\cdot y=f^1 &\mbox{ in }\Omega,\\\ns -\displaystyle\sum_{p,
q=1}^{m}\left[(a_{p q}^{21}y^1_{x_p})_{x_q}+(a_{p
q}^{22}y^2_{x_p})_{x_q}+\cdots+(a_{p
q}^{2n}y^n_{x_p})_{x_q}\right]+\displaystyle\sum_{i=1}^nC^{2i}\cdot\nabla
y^i+D^2\cdot y=f^2 &\mbox{ in }\Omega,\\\ns
\quad\quad\quad\quad\quad\quad\quad\quad\quad
\quad\quad\quad\quad\quad\quad\quad\vdots\\\ns
-\displaystyle\sum_{p, q=1}^{m}\left[(a_{p
q}^{n1}y^1_{x_p})_{x_q}+(a_{p q}^{n2}y^2_{x_p})_{x_q}+\cdots+(a_{p
q}^{nn}y^n_{x_p})_{x_q}\right]+\displaystyle\sum_{i=1}^nC^{ni}\cdot\nabla
y^i+D^n\cdot y=f^n &\mbox{ in }\Omega,\\\ns
y^1=g^1,\quad
y^2=g^2,\quad \cdots,\quad y^n=g^n  &\mbox{ on }\Gamma.
\end{array}
\right.
\end{eqnarray}
In both (\ref{1}) and (\ref{61}), $y=(y^1, \cdots, y^n)^\top$ is
unknown, while $a^{ij}$, $a^{i j}_{p q}$, $C^{i j}$, $D^i$, $f^i$
and $g^i$ ($i, j=1, \cdots, n$; $p, q=1, \cdots, m$) are suitable
given functions (See the next section for the assumptions on these
functions). The main purpose of this paper is to study the weak
maximum principle, or the boundedness of weak solutions, for systems
(\ref{1}) and (\ref{61}) with suitable measurable principal part
coefficients.

It is well-known that the weak maximum principle is one of the basic
issues in the theory of partial differential equations and it plays
an essential role in the study of many other problems. For example,
a central problem in the calculus of variations is the regularity of
stationary points for functionals of the type
 $$J(u) =
\int_\Omega
 F(x, u(x),\nabla u(x)) dx,
 $$
where $u(x) = (u^1(x),u^2(x),\cdots,u^n(x))^\top$ is a vector-valued
function defined on $\Omega$, and $F(\cdot, \cdot, \cdot)$ is a
suitable function defined on $\Omega\times \dbR^n\times \dbR^{mn}$.
Research in this area has been stimulated by D.~ Hilbert's Problem
19, which can be reduced to the regularity of weak solutions to
nonlinear elliptic equations or systems. This problem was
successfully solved by C. Morrey (\cite{Morrey}) in two dimensions
and the general case with $n = 1$ was finally solved by E. De Giorgi
(\cite{de}) and J. Nash (\cite{n}), and refined by J. Moser
(\cite{m}). We refer to \cite{b,BB,g,gt,l} and the references cited
therein for more details in this respect. A fundamental step of the
De Giorgi-Nash-Moser approach in solving the scalar Hilbert's 19th
problem (i.e., $n = 1$) is to establish the weak maximum principle
for single linear elliptic equations.

In many physical and geometrical applications, $u$ may be a vector
function, and therefore, the corresponding Euler-Lagrange equation
is a system. Naturally, one expects to extend the De
Giorgi-Nash-Moser approach to the case of systems. However, in 1968,
E. De Giorgi (\cite{d}) gave a surprising counterexample of an
unbounded solution to a second order linear elliptic system with
bounded coefficients. This means that the weak maximum principle
fails for general second order linear elliptic systems, and
therefore, the De Giorgi-Nash-Moser estimates are no longer valid
for general elliptic systems.

In order to establish the weak maximum principle for elliptic
differential systems, as a consequence of the above mentioned De
Giorgi's counterexample, one has to impose some restrictions on the
structure of the system. There exist a few works in this direction.
In \cite{l}, a weak maximum principle was proved for a class of
special elliptic systems with variable coefficients, in which the
principal operator in each equation takes the same form, and it is
acting only on one component of the solution vector. In
\cite{cc,Can,le}, some weak maximum principles were discussed in the
frame of Campanato's space for linear or quasilinear elliptic
systems under some additional conditions, say, $2\leq m\leq 4$ in
\cite{cc}, the coefficients matrix being constant in \cite{Can}, and
a dispersion assumption on the eigenvalues of the principal part
coefficients matrix in \cite{le} (and hence the system is
almost-diagonal in high space dimensions).

In this paper, we choose the usual Sobolev space as the working
space and derive weak maximum principles for two classes of strongly
coupled elliptic systems with different principal parts, in the
spirit of the classical framework for single equations. We emphasize
that our systems are strongly coupled, i.e., the (second order)
terms of the principal parts are coupled each other. Therefore, when
establishing the desired {\it a prior} estimate, it is necessary to
get rid of some undesired terms generated by different principal
operators and/or different solution components appeared in the same
equation. This goal is achieved by choosing delicately suitable
weighted test functions. As far as we know, this is the first result
on the weak maximum principle (in the classical sense) for strongly
coupled elliptic systems.

The rest of this paper is organized as follows. Section 2 is devoted
to stating the main results in this work. In Section 3, we collect
some preliminary results which will be useful later. Sections 4 and
5 are addressed to the proof of the main results, i.e., the
boundedness of weak solutions to systems (\ref{1}) and (\ref{61}),
respectively. Finally, in Section 6, we give an example in which the
assumptions for proving the boundedness of the weak solution to
system (\ref{61}) are satisfied.

\section{Statement of the main results}

To begin with, we introduce some assumptions. Suppose that, for $i,
j=1, \cdots, n$,
 \bel{2e1}
 a^{i j} \in L^\infty(\Omega)
 \ee
and
 \bel{1co2}
 \left\{
 \ba{ll}
\ds C^{i j}(\cdot)\in L^\theta(\Omega; \dbR^m)\ \hb{ and }\
D^i(\cdot)\in L^{\frac{\theta}{2}}(\Omega; \dbR^n)\q \hb{for some
}\theta>m,
 \\
 \ns
 \ds f=(f^1, \cdots, f^n)^\top\in H^{-1}(\Omega;\dbR^n),
 \q g=(g^1, \cdots, g^n)^\top\in H^1(\Omega;\dbR^n),
 \ea
\right.
 \ee
and for $i, j=1, \cdots, n$ and $p, q=1, \cdots, m$,
 \bel{2co2}
 a^{i j}_{p q} \in L^\infty(\Omega),\qquad  a_{p q}^{i j}=a_{q p}^{i j}.
 \ee
Moreover, we assume that, for some positive constant $\rho$,
\begin{equation}\label{2}
 \ba{ll}
\displaystyle\sum_{i, j=1}^{n}\displaystyle\sum_{p=1}^{m}a^{i
j}\xi_p^i\xi_p^j\geq \rho|\xi|^2,\\
 \ns
 \ds\qq\qq\qq\qq\qq\forall\;(x, \xi)=(x, \xi^1_1, \cdots, \xi^1_m, \cdots, \xi^n_1,
\cdots, \xi^n_m)\in \Omega\times \dbR^{nm},
 \ea
\end{equation}
and
\begin{equation}\label{63}
 \ba{ll}
\displaystyle\sum_{i, j=1}^{n}\displaystyle\sum_{p, q=1}^{m}a_{p q}^{i j}(x)\xi^i_p \xi^j_q \geq \rho|\xi|^2,\\
 \ns
 \ds\qq\qq\qq\qq\qq\forall\;(x, \xi)=(x, \xi^1_1, \cdots, \xi^1_m, \cdots, \xi^n_1,
\cdots, \xi^n_m)\in \Omega\times \dbR^{nm}.
 \ea
\end{equation}

Conditions (\ref{2}) and (\ref{63}) mean that both systems (\ref{1})
and (\ref{61}) are elliptic (see \cite[Section 1 of Chapter 8]{c}).
Clearly, system (\ref{1}) is a special case of system (\ref{61}).
The weak solution to system (\ref{61}) is understood in the
following sense:

\begin{definition}\label{80}
We call $y=(y^1, \cdots, y^n)^\top\in H^1(\Omega; \dbR^n)$ to be a
weak solution to system $(\ref{61})$ if for any $\varphi=(\varphi^1,
\cdots, \varphi^n)^\top \in H^1_0(\Omega; \dbR^n)$,
\begin{eqnarray*}
&&\displaystyle\sum_{i, j=1}^{n}\displaystyle\sum_{p,
q=1}^{m}\displaystyle\int_\Omega a^{i j}_{p q}(x)
y^j_{x_p}\varphi^i_{x_q}
dx+\displaystyle\int_\Omega\displaystyle\sum_{i=1}^{n}
\Big[\displaystyle\sum_{j=1}^{n} C^{ij}(x)\cdot\nabla
y^j\varphi^i+D^i(x)\cdot y \varphi^i\Big]dx\\\ns
&&=\langle f,
\varphi\rangle_{H^{-1}(\Omega; \dbR^n), H^1_0(\Omega; \dbR^n)},
\end{eqnarray*}
and $y^i-g^i\in H^1_0(\Omega),\ i=1, \cdots, n$.
\end{definition}

Similar to the proof of \cite[Theorem 2.3 in Chapter 1]{c}), it is
easy to show the following well-posedness result for system
(\ref{61}).

\bl\label{2l1}
Let conditions $(\ref{1co2})$, $(\ref{2co2})$ and $(\ref{63})$ be
fulfilled. Then, there exists a constant $\nu_0=\nu_0(n,m,\theta)>0$
such that system $(\ref{61})$ admits a unique weak solution $y\in
H^1(\Omega; \dbR^n)$ whenever the following inequality
 \begin{equation}\label{65}
 \ba{ll}
\displaystyle\sum_{i=1}^n \left(D^i(x)\cdot \mu\right) \mu^i\geq
\nu\rho^{\frac{m+\theta}{m-\theta}}\Big[\displaystyle\sum_{i, j=1}^{n} |C^{ij}|_{L^\theta(\Omega; \dbR^m)}\Big]^{\frac{2\theta}{\theta-m}}|\mu|^2,\\
 \ns
 \ds\qq\qq\qq\qq\forall\; (x, \mu)=(x, \mu^1, \mu^2, \cdots, \mu^n)\in
\Omega\times\dbR^n
 \ea
\end{equation}
holds for $\nu\ge \nu_0$. Moreover,
$$
|y|_{H^1(\Omega; \dbR^n)}\leq C\Big(n, m, \Omega, \rho, |a_{p q}^{i
j}|_{L^\infty(\Omega)}, |C^{i j}|_{L^\theta(\Omega; \dbR^m)},
|D^i|_{L^{\frac{\theta}{2}}(\Omega;
\dbR^n)}\Big)(|f|_{H^{-1}(\Omega; \dbR^n)}+|g|_{H^1(\Omega;
\dbR^n)}).
$$
\el

The proof of Lemma \ref{2l1} is standard and therefore we omit the
details.

\br
Clearly, if $C^{i j}(\cdot)\equiv 0$ for $i,j=1,\cdots,n$, then
condition (\ref{65}) is satisfied whenever the function matrix
$\big(D^1(x),\cdots, D^n(x)\big)$ is semi-positive definite.
\er

Next, we put
$$
A= \left(
\begin{array}{lccl}
a^{11} & a^{21} & \cdots & a^{n1} \\
a^{12} & a^{22} & \cdots & a^{n2} \\
\vdots & \vdots & \vdots & \vdots \\
a^{1n} & a^{2n} & \cdots & a^{nn}
\end{array}
\right), \qquad \ \ \ B= \left(
\begin{array}{lccl}
a^{22} & a^{32} & \cdots & a^{n2} \\
a^{23} & a^{33} & \cdots & a^{n3} \\
\vdots & \vdots & \vdots & \vdots \\
a^{2n} & a^{3n} & \cdots & a^{nn}
\end{array}
\right).
$$
Also, denote by $B^{ij}$ $(i, j=1, \cdots, n)$ the cofactor of $A$
with respect to $a^{i j}$ and by $\det A$ the determinant of matrix
$A$. It is easy to see that $B^{1 1}=B$. Moreover, under condition
$(\ref{2})$, it is easy to show that $\det B\neq 0$.

The first main result in this paper is the following boundedness of
weak solutions to system (\ref{1}).

\begin{theorem}\label{5}
Suppose that conditions $(\ref{2e1})$, $(\ref{1co2})$ and
$(\ref{2})$ are fulfilled, inequality (\ref{65}) holds for $\nu\ge
\nu_0$ (Recall Lemma \ref{2l1} for $\nu_0$), $f\in
L^\frac{\theta}{2}(\Omega; \dbR^n)$ and
 \bel{kkk1}
 \frac{\det B^{i j}}{\det B} \in W^{1,
\infty}(\Omega),\qq i, j=1, \cdots, n.
  \ee
Then the weak solution $y\in H^1(\Omega; \dbR^n)$ to system
$(\ref{1})$ satisfies
\begin{eqnarray*}
&&\esssup\limits_\Omega |y|\leq C\left(m,\ n,\ \theta,\ \Omega, \
\rho,\ |a^{i j}|_{L^\infty(\Omega)},\ |C^{i j}|_{L^\theta(\Omega;
\dbR^m)},\ |D^i|_{L^\frac{\theta}{2}(\Omega; \dbR^n)},\
\left|\displaystyle\frac{\det B^{i j}}{\det B}\right|_{W^{1, \infty}(\Omega)},\right.\\
&&\left.\quad\quad\quad\quad\quad\quad\quad |g|_{H^1(\Omega;
\dbR^n)},\ |f|_{L^\frac{\theta}{2}(\Omega; \dbR^n)}, \
\esssup\limits_{\Gamma}|y|\right).
\end{eqnarray*}
\end{theorem}

The proof of Theorem \ref{5} will be given in Section \ref{s4}.

\br
We conjecture that assumption (\ref{kkk1}) in Theorem \ref{5} is a
technical condition, and therefore it is not really necessary.
However, we do not know how to drop this assumption at this moment.
\er

Since almost all of the natural materials are isotropic, Theorem
\ref{5} suffices for most of physical applications. Nevertheless,
from the mathematical point of view, it would be quite interesting
to extend Theorem \ref{5} to more general anisotropic systems such
as $(\ref{61})$, in which the scalar functions $a^{ij}$
($i,j=1,2,\cdots,n$) appeared in system $(\ref{1})$ are replaced by
the $\mathbb{R}^{m\times m}$ matrix-valued functions
$\big(a^{ij}_{pq}\big)_{1\le p,q\le m}$. Note however that, by the
above mentioned De Giorgi's counterexample (\cite{d}), this seems to
be highly nontrivial in the general case. In the rest of this
section, we shall extend Theorem \ref{5} to system $(\ref{61})$
under some technical assumptions.

In order to treat system $(\ref{61})$, we put
$$
M_{p q}= \left(
\begin{array}{lccl}
a_{p q}^{11} & a_{p q}^{21} & \cdots & a_{p q}^{n1} \\[2mm]
a_{p q}^{12} & a_{p q}^{22} & \cdots & a_{p q}^{n2} \\[2mm]
\vdots & \vdots & \vdots & \vdots \\
a_{p q}^{1n} & a_{p q}^{2n} & \cdots & a_{p q}^{nn}
\end{array}
\right), \qquad \ \ \  L_{p q}=\mbox{det }M_{p q}, \qquad \  (p,
q=1, \cdots, m).
$$
We assume that
 \bel{2e7}
L_{p q}\neq 0,\qq \forall\;p, q=1, \cdots, m.
 \ee
Also, we denote by $v^{i j}_{p q}$ $(i, j=1, \cdots, n; \ \ p, q=1,
\cdots, m)$ the cofactor of $M_{p q}$ with respect to $a^{i j}_{p
q}$.

Further, let us introduce the following assumption:

\ms\it

\no{\bf (H)} There exist functions $f_{p q},\ h^{i j}\in W^{1,
\infty}(\Omega)$ $(i, j=1, \cdots, n;\ \ p, q=1, \cdots, m)$ such
that

\begin{enumerate}

\item[{\bf 1)}] $h^{11}\equiv 1$, $h^{i j}=h^{j i}$, and
the following matrix is uniformly positive definite:
$$V:=\left(
\begin{array}{lccl}
1 & h^{1 2} & \cdots & h^{1 n} \\[2mm]
h^{2 1} &  h^{2 2} & \cdots & h^{2 n} \\[2mm]
\vdots & \vdots & \vdots & \vdots \\
h^{n 1} & h^{n 2} & \cdots & h^{n n}
\end{array}
\right),
$$
i.e., $V\geq\ \rho_1 I_{n\times n}$ for some positive number
$\rho_1$;

\item[{\bf 2)}] The function
$E^{i j}:=\displaystyle\frac{f_{p q}}{L_{p q}}
\displaystyle\sum_{l=1}^{n} h^{l j}v^{l i}_{p q}$ is independent
of $p$ and $q$, and $E^{i j}\in W^{1, \infty}(\Omega)$ for any $i,
j=1, \cdots, n$;

\item[{\bf 3)}] The following matrix is uniformly positive definite:
$$F:=\left(
\begin{array}{lccl}
F_{1 1} & F_{1 2} & \cdots & F_{1 m} \\[2mm]
F_{2 1} &  F_{2 2} & \cdots & F_{2 m} \\[2mm]
\vdots & \vdots & \vdots & \vdots \\
F_{m 1} & F_{m 2} & \cdots & F_{m m}
\end{array}
\right),
$$
i.e., $F\geq\ \rho_2 I_{m\times m}$ for some positive number
$\rho_2$, where $F_{p q}:=\displaystyle\sum_{l=1}^{n}a^{l 1}_{p
q}E^{l 1}$ for any $p, q=1, \cdots, m$;

\item[{\bf 4)}] The following matrix is uniformly positive definite:
$$M:=\left(
\begin{array}{lccl}
F & h^{1 2}F & \cdots & h^{1 n}F \\ [2mm]
h^{2 1}F &  h^{2 2}F & \cdots & h^{2 n}F \\[2mm]
\vdots & \vdots & \vdots & \vdots \\
h^{n 1}F & h^{n 2}F & \cdots & h^{n n}F
\end{array}
\right)_{nm\times nm},
$$
i.e., $M\geq\ \rho_3 I_{nm\times nm}$ for some positive number
$\rho_3$.

\end{enumerate}

\rm\ms

Now, we can state our another main result as the following
boundedness result for weak solutions to system (\ref{61}).

\begin{theorem}\label{64}
Suppose that conditions  $(\ref{2co2})$, $(\ref{63})$ and
$(\ref{2e7})$ are fulfilled, $C^{i j}(\cdot)\in L^\infty(\Omega;
\dbR^m)$, $ D^i(\cdot)\in L^{\infty}(\Omega; \dbR^n)$ and inequality
(\ref{65}) holds for $\nu\ge \nu_0$ ((Recall Lemma \ref{2l1} for
$\nu_0$)), $f\in L^\infty(\Omega; \dbR^n)$, and assumption {\rm \bf
(H)} holds. Then the weak solution $y\in H^1(\Omega; \dbR^n)$ to
system $(\ref{61})$ satisfies
\begin{eqnarray*}
&&\esssup\limits_\Omega |y|\leq C\left(m,\ n,\ \Omega, \ \rho,\
\rho_1,\ \rho_2,\ \rho_3,\ |a^{i j}_{p q}|_{L^\infty(\Omega)},\
|C^{i j}|_{L^\infty(\Omega; \dbR^m)},\ |D^i|_{L^\infty(\Omega;
\dbR^n)},\
\right.\\
&&\left.\quad\quad\quad\quad\quad\quad\quad \left|E^{i
j}\right|_{W^{1, \infty}(\Omega)},\ |h^{i j}|_{W^{1,
\infty}(\Omega)}, |g|_{H^1(\Omega; \dbR^n)},\ |f|_{L^\infty(\Omega;
\dbR^n)}, \ \esssup\limits_{\Gamma}|y|\right).
\end{eqnarray*}
\end{theorem}

The proof of Theorem \ref{64} will be given in Section \ref{s5}.
Also, in Section \ref{s6}, we shall give an illustrative example, in
which all of the assumptions in Theorem \ref{64} are satisfied.

\begin{remark}
It is well-known that one of the classical topics in partial
differential equations is the strong maximum principle for elliptic
differential equations, which has many applications
(\cite{Fraenkel,Nir,PW,PS} and so on). However, the existing results
on strong maximum principle are mainly focusing on single elliptic
equations, although one can find some works on weakly coupled
elliptic systems (\cite{b,JM,Sirakov}) and the references therein.
It would be quite interesting to establish a strong maximum
principle for system $(\ref{1})$ or even for system $(\ref{61})$,
but this remains to be done and it seems to be far from easy.
\end{remark}

\section{Some preliminaries}

In this section, we collect some known preliminary results which
will be useful later.

The first one is the following interpolation result.

\begin{lemma}\label{31} {\rm (\cite[Theorem $2.1$ in Chapter $2$]{l})} For any $u \in
W^{1,t}_0(\Omega),$ $t\geq 1$ and $\tau\geq 1$, it holds that
$$|u|_{L^{p^*}(\Omega)}\leq \beta |\nabla
u|^\alpha_{L^t(\Omega)}|u|^{1-\alpha}_{L^\tau(\Omega)},$$ where
$\alpha=\left(\frac{1}{\tau}-\frac{1}{p^*}\right)\left(\frac{1}{\tau}-\frac{1}{t^*}\right)^{-1}$,
$t^*=\frac{t m}{m-t}$, and $\beta$ is a constant depending only on
$m$, $t$, $p^*$, $\tau$ and $\alpha$. Moreover, if $t<m$, $p^*$ can
be any number between $\tau$ and $t^*$; if $t\geq m$, $p^*$ can be
any number larger than $\tau$.
\end{lemma}

For any Lebesgue measurable function $u$ defined on $\Omega$, we put
$A_k=\{x\in \Omega; u(x)> k\}$ and denote by $|A_k|$ the Lebesgue
measure of set $A_k$. The next lemma is quite useful in deriving the
supremum of function $u$.

\begin{lemma}\label{32} {\rm (\cite[Theorem $5.1$ in Chapter $2$]{l})}
Suppose that $u\in W^{1,m_0}(\Omega)\cap L^{q_0}(\Omega)$ for some
$m_0\in [1, m]$ and some $q_0\geq 1$. If for any fixed $k\geq
\esssup\limits_{\Gamma}u$, function $u$ satisfies the following
inequality:
$$\displaystyle\int_{A_k}|\nabla u|^{m_0} dx \leq \gamma
\left[\displaystyle\int_{A_k}(u-k)^{l_0}
dx\right]^{\frac{m_0}{l_0}}+\gamma
k^\sigma|A_k|^{1-\frac{m_0}{m}+\varepsilon_0},
 $$
where $\gamma$, $l_0$, $\sigma$ and $\varepsilon_0$ are positive
constants satisfying $l_0<\frac{mm_0}{m-m_0}$ and $m_0\leq
\sigma<\varepsilon_0 q_0+m_0$, then
$$\esssup\limits_{\Omega} u\leq C^*(\Omega,\ m_0,\ q_0,\
\gamma,\ l_0,\ \sigma,\ \varepsilon_0,\ \esssup\limits_{\Gamma}u,\
|u|_{L^{q_0}(\Omega)}).
 $$
Moreover, when $\sigma=m_0$, $|u|_{L^{q_0}(\Omega)}$ appeared in
$C^*$ can be replaced by $|u|_{L^1(\Omega)}$.
\end{lemma}

The last lemma is a result on comparison of the determinants between
a matrix and its symmetrizing matrix.

\begin{lemma}\label{33} {\rm (\cite[Theorem $3.7.1$] {w})} For a real
matrix $E$, if $H(E)=\displaystyle\frac{E+E^\top}{2}$ is positive
definite, then
$$\det H(E)\leq \det E.$$
\end{lemma}

\section{Proof of Theorem \ref{5}}\label{s4}
The goal of this section is to prove our first main result, i.e.,
Theorem \ref{5}.\\

\noindent {\bf Proof of Theorem \ref{5}.} First, for the weak
solution $y=(y^1, y^2, \cdots, y^n)^\top$ to system $(\ref{1})$, any
fixed $k\geq \esssup\limits_{\Gamma}|y|^2$ and $r>0$, put
$$\phi_r(x)=\mbox{min}\{(|y(x)|^2-k)_+, r\},$$ where $s_+=\max\{s,
0\}$ (for any $s\in\dbR$). We choose $\varphi=(\varphi^1, \varphi^2,
\cdots, \varphi^n)^\top$ as a test function, where
$\varphi^i=(y^1T^{1 i}+y^2T^{2 i}+\cdots+y^nT^{n i})\phi_r$ and
$T^{i j}$ $(i, j=1,2, \cdots, n)$ are suitable functions to be
specified later. Obviously, $\varphi\in H^1_0(\Omega; \dbR^n)$. By
Definition \ref{80}, it follows that
\begin{eqnarray*} &&\displaystyle\sum_{i=1}^n\displaystyle\int_\Omega(a^{i1}\nabla y^1
+a^{i2}\nabla y^2+\cdots+a^{in}\nabla y^n)\cdot\nabla
\left[(y^1T^{1 i}+y^2T^{2 i}+\cdots+y^nT^{n i})\phi_r\right]dx\\
&&\q+\displaystyle\sum_{i=1}^n\displaystyle\int_\Omega\Big(
\displaystyle\sum_{j=1}^n C^{i j}\cdot \nabla y^j+D^i\cdot
y\Big)(y^1T^{1 i}+y^2T^{2 i}+\cdots+y^nT^{n i})\phi_rdx\\
 &&=\displaystyle\sum_{i=1}^n\displaystyle\int_\Omega
f^i(y^1T^{1 i}+y^2T^{2 i}+\cdots+y^nT^{n i})\phi_rdx.
\end{eqnarray*}
This implies that, for any $\e_1>0$,
 \be\label{7}
\begin{array}{ll}
\displaystyle\sum_{i=1}^n\displaystyle\int_\Omega \Big\{ \left(a^{i
1}T^{1 i}|\nabla y^1|^2+a^{i 2}T^{2 i}|\nabla y^2|^2+\cdots +a^{i
n}T^{n i}|\nabla y^n|^2\right)\phi_r\\\ns
\quad\quad+\displaystyle\sum_{l, j\in \{1, 2, \cdots, n\}, \,l\not=
j} \left(a^{i j}T^{l i}\nabla y^j\cdot\nabla y^l \phi_r+a^{i j} T^{l
i} y^l\nabla y^j\cdot\nabla\phi_r\right)\\\ns
\quad\quad+\displaystyle\frac{1}{2}\left[ a^{i 1}T^{1
i}\nabla(y^1)^2+ a^{i 2}T^{2 i}\nabla(y^2)^2+\cdots+a^{i n}T^{n
i}\nabla(y^n)^2\right]\cdot\nabla\phi_r\Big\}dx\\\ns
=\displaystyle\sum_{i=1}^n\displaystyle\int_\Omega\Big[ f^i(y^1T^{1
i}+y^2T^{2 i}+\cdots+y^nT^{n i})\phi_r-\displaystyle\sum_{l,
j=1}^{n}a^{i l}y^j\nabla y^l\cdot\nabla T^{j i} \phi_r\\\ns
\quad\quad-\Big( \displaystyle\sum_{j=1}^n C^{i j}\cdot \nabla
y^j+D^i\cdot y\Big)(y^1T^{1 i}+y^2T^{2 i}+\cdots+y^nT^{n
i})\phi_r\Big]dx\\\ns
\leq C\displaystyle\int_\Omega|f||y|\phi_r
dx+\e_1\displaystyle\int_{\Omega}|\nabla y|^2\phi_r
dx+C\displaystyle\int_{\Omega}\left[\e_1^{-1}\left(\sum_{i,j=1}^n|C^{i
j}|^2+1\right)+\sum_{i=1}^n|D^i|\right]|y|^2\phi_r dx.
\end{array}
\ee
Here and hereafter $C$ denotes a generic constant (which may be
different from one place to another),  depending only on $n$, $m$,
$\rho$, $|a^{i j}|_{L^\infty(\Omega)}$ and $|T^{i j}|_{W^{1,
\infty}(\Omega)}$ $(i, j=1, 2, \cdots, n)$. Also, it is clear that
\begin{eqnarray}\label{8}
\begin{array}{rl}
&\left[ a^{i 1}T^{1 i}\nabla(y^1)^2+ a^{i 2}T^{2
i}\nabla(y^2)^2+\cdots+a^{i
n}T^{n i}\nabla(y^n)^2\right]\cdot\nabla\phi_r\\[3mm]
&=a^{i 1}T^{1 i}\nabla|y|^2\cdot\nabla\phi_r+ \left[(a^{i 2}T^{2
i}-a^{i 1}T^{1 i})\nabla(y^2)^2+\cdots+(a^{i n}T^{n i}-a^{i 1}T^{1
i})\nabla(y^n)^2\right]\cdot\nabla\phi_r.
\end{array}
\end{eqnarray}

Next, we choose suitable functions $T^{i j}$ $(i, j=1, 2, \cdots,
n)$ such that $T^{i j}\in W^{1, \infty}(\Omega)$ and some undesired
terms in (\ref{7}) vanish. For this purpose, we consider the
following linear system:
\begin{eqnarray}\label{22}
\left\{
\begin{array}{rl}
&\displaystyle\sum_{i=1}^{n}a^{i j}T^{l i}=0,\qq \forall\;j,\ l=1,
2, \cdots, n\hb{ with }j\neq l,\\\ns
&\displaystyle\sum_{i=1}^{n}(a^{i 2}T^{2 i}-a^{i
1}T^{1 i})=0,\\
&\quad\quad\quad\quad\quad\vdots\\
&\displaystyle\sum_{i=1}^{n}(a^{i n}T^{n i}-a^{i 1}T^{1 i})=0.
\end{array}
\right.
\end{eqnarray}
By Lemma \ref{33} and (\ref{2}), it follows that $\det A\geq \rho^n$
and $\det B\geq \rho^{n-1}$. One can easily check that the following
is a solution to system (\ref{22}):
\begin{eqnarray}\label{34}
 \left\{
 \begin{array}{rl}
&T^{1 1}=1, \ \ \ \ \ T^{1 i}=(-1)^{1+i}\displaystyle\frac{\det B^{1
i}}{\det B}\ \  (i=2, \cdots, n),\\\ns
 &T^{j i}=(-1)^{i+j}
\displaystyle\sum_{l=1}^{n}a^{l 1}T^{1 l} \displaystyle\frac{\det
B^{j i}}{\det A}\ \  (j=2, 3, \cdots, n,\ \ i=1, 2, \cdots, n).
\end{array}\right.
\end{eqnarray}
Also, it is not difficult to check that
$$
\displaystyle\sum_{l=1}^{n}a^{l 1}T^{1 l}=\displaystyle\frac{\det
A}{\det B}.
$$
From this fact and noting (\ref{34}), we see that $T^{j
i}=(-1)^{i+j}\displaystyle\frac{\det B^{j i}}{\det B}$ $(i, j=1, 2,
\cdots, n)$. Moreover, there exists a constant $\rho^*>0$, depending
only on $n, \rho$ and $|a^{i j}|_{L^\infty(\Omega)}$ $(i,
j=2,\cdots, n)$, such that
\begin{equation}\label{35}
\displaystyle\sum_{l=1}^{n}a^{l 1}T^{1 l}=\cdots=
\displaystyle\sum_{l=1}^{n}a^{l n}T^{n l}\geq
\displaystyle\frac{\rho^n}{\det B}\geq  \rho^*.
\end{equation}

Combining (\ref{8}), (\ref{22}) and (\ref{35}) with (\ref{7}), we
arrive at
 \bel{14} \displaystyle\int_\Omega\left(|\nabla y|^2\phi_r
+|\nabla\phi_r|^2\right)dx \leq
C\left[\displaystyle\int_\Omega|f||y|\phi_r
dx+\displaystyle\int_{\Omega}\left(\sum_{i,j=1}^n|C^{i
j}|^2+\sum_{i=1}^n|D^i|+1\right)|y|^2\phi_r dx\right].
 \ee

In the sequel, we estimate the terms in the right side of (\ref{14})
and prove that the left side of this inequality is uniformly bounded
with respect to $r>0$. First of all, by H\"older's inequality and
Lemma \ref{31}, we see that
\begin{eqnarray}\label{17}
\begin{array}{rl}
&\displaystyle\int_\Omega|f||y|\phi_r dx \leq
\displaystyle\int_\Omega |f|(|y|^2-k)^{\frac{1}{2}}\phi_r
dx+k^{\frac{1}{2}}\displaystyle\int_\Omega |f|\phi_r dx\\[3mm]
&\leq  \left|f\right|_{L^{\frac{4\theta}{\theta+6}}(\Omega;
\dbR^n)}\left|(|y|^2-k)^
{\frac{1}{2}}\phi_r\right|_{L^{\frac{4\theta}{3\theta-6}}(\Omega)}
+C_1\left|f\right|_{L^{\frac{\theta}{2}}(\Omega;
\dbR^n)}\left|\phi_r\right|_{L^{\frac{\theta}{\theta-2}}(\Omega)}\\[3mm]
&\leq  \left|f\right|_{L^{\frac{4\theta}{\theta+6}}(\Omega;
\dbR^n)}\left|
(|y|^2-k)\phi_r\right|^{\frac{3}{4}}_{L^{\frac{\theta}{\theta-2}}(\Omega)}
+C_1\left|f\right|_{L^{\frac{\theta}{2}}(\Omega;
\dbR^n)}|y|^2_{H^1(\Omega; \dbR^n)}\\[3mm]
&\leq \left|f\right|_{L^{\frac{4\theta}{\theta+6}}(\Omega;
\dbR^n)}\Big[\left|
(|y|^2-k)\phi_r\right|_{L^{\frac{\theta}{\theta-2}}(\Omega)}+1\Big]
+C_1\left|f\right|_{L^{\frac{\theta}{2}}(\Omega;
\dbR^n)}|y|^2_{H^1(\Omega; \dbR^n)},
\end{array}
\end{eqnarray}
here and hereafter $C_1$ denotes a constant  depending only on $k$,
$\theta$, $n$, $m$ and $\Omega$. Put $\ds L=\sum_{i,j=1}^n|C^{i
j}|_{L^\theta(\O;\dbR^m)}^2+\sum_{i=1}^n|D^i|_{L^{\frac
{\theta}{2}}(\O;\dbR^n)}+1$. By Lemma \ref{31}, we find that
\begin{eqnarray}\label{45}
\begin{array}{rl}
&\displaystyle \int_{\Omega}\left(\sum_{i,j=1}^n|C^{i
j}|^2+\sum_{i=1}^n|D^i|+1\right)|y|^2\phi_r dx\\\ns
 &\leq C_1 L
\left[\displaystyle\int_\Omega
(|y|^2\phi_r)^{\frac{\theta}{\theta-2}}dx\right]^{\frac{\theta-2}{\theta}}\leq
C_1 L\left|
(|y|^2-k)\phi_r\right|_{L^{\frac{\theta}{\theta-2}}(\Omega)}+ C_1 L
\left|\phi_r\right|_{L^{\frac{\theta}{\theta-2}}(\Omega)}\\\ns &\leq
C_1 L\left|
(|y|^2-k)\phi_r\right|_{L^{\frac{\theta}{\theta-2}}(\Omega)}+
C_1L|y|^2_{H^1(\Omega; \dbR^n)}.
\end{array}
\end{eqnarray}
On the other hand, put $u_*=\sqrt{(|y|^2-k)\phi_r}$. It follows that
\begin{eqnarray}\label{19}
\begin{array}{rl}
&\displaystyle\int_\Omega|\nabla u_*|^2dx=\int_{|y|^2> k}|\nabla
u_*|^2dx =\int_{|y|^2> k}\left|\frac{2\phi_ry\nabla y
+(|y|^2-k)\nabla\phi_r}{2\sqrt{(|y|^2-k)\phi_r}}\right|^2dx
\\\ns
&\ds\leq 2\int_{|y|^2> k}\left|\frac{\phi_ry\nabla y
}{\sqrt{(|y|^2-k)\phi_r}}\right|^2dx
+\displaystyle\frac{1}{2}\int_{k+r\ge |y|^2>
k}\left|\frac{(|y|^2-k)\nabla\phi_r}{\sqrt{(|y|^2-k)\phi_r}}\right|^2dx\\\ns
&\leq 2\displaystyle\int_{|y|^2> k} (|y|^2-k)^{-1}\phi_r|y|^2|\nabla
y|^2 dx+2\displaystyle\int_{k+r\ge|y|^2>
k}(|y|^2-k)\phi_r^{-1}|y|^2|\nabla y|^2 dx.
\end{array}
\end{eqnarray}
Noting $\phi_r\le |y|^2-k$, we see that
 \bel{ops1}
 \ba{ll}
 \displaystyle\int_{|y|^2> k} (|y|^2-k)^{-1}\phi_r|y|^2|\nabla
y|^2 dx=\int_{|y|^2> k} \phi_r|\nabla y|^2 dx+k\int_{|y|^2> k}
(|y|^2-k)^{-1}\phi_r|\nabla y|^2 dx\\\ns
 \ns
 \ds\le \int_\O \phi_r|\nabla y|^2 dx+k\int_{|y|^2> k}
|\nabla y|^2 dx\le \int_\O \phi_r|\nabla y|^2 dx+k\int_\O |\nabla
y|^2 dx.
 \ea
 \ee
Noting that $\phi_r= |y|^2-k$ whenever $k+r\ge |y|^2$,  it is clear
that
 \bel{ops2}
 \ba{ll}\ds
 \int_{k+r\ge|y|^2>
k}(|y|^2-k)\phi_r^{-1}|y|^2|\nabla y|^2 dx=\int_{k+r\ge|y|^2>
k}|y|^2|\nabla y|^2 dx\\\ns
 \ns
 \ds=\int_{k+r\ge|y|^2>
k}(|y|^2-k)|\nabla y|^2 dx+k\int_{k+r\ge|y|^2> k}|\nabla y|^2
dx\\\ns
 \ns
 \ds=\int_{k+r\ge|y|^2>
k}\phi_r|\nabla y|^2 dx+k\int_{k+r\ge|y|^2> k}|\nabla y|^2 dx
 \le
\int_\O \phi_r|\nabla y|^2 dx+k\int_\O |\nabla y|^2 dx.
 \ea
 \ee
Therefore, by (\ref{19})--(\ref{ops2}), we conclude that
 \bel{1oo9}
\int_\Omega|\nabla u_*|^2dx\leq 4\displaystyle\int_\Omega|\nabla
y|^2\phi_r dx+C_1\displaystyle\int_\Omega|\nabla y|^2dx. \ee

By (\ref{1oo9}) and Lemma \ref{31}, for any $0<\varepsilon_2<1$, we
end up with
\begin{eqnarray}\label{20}
\begin{array}{rl}
&\left|
(|y|^2-k)\phi_r\right|_{L^{\frac{\theta}{\theta-2}}(\Omega)}=\left(\displaystyle\int_\Omega
u_*
^{\frac{2\theta}{\theta-2}}dx\right)^{\frac{\theta-2}{\theta}}\leq
\varepsilon_2\displaystyle\int_\Omega|\nabla
u_*|^2dx+C_1\varepsilon_2^{-1}\left(\displaystyle\int_\Omega|u_*|dx\right)^2\\\ns
&\leq 4\varepsilon_2\displaystyle\int_\Omega|\nabla y|^2\phi_r
dx+C_1\varepsilon_2|y|^2_{H^1(\Omega;
\dbR^n)}+C_1\varepsilon_2^{-1}\left[\displaystyle\int_\Omega
(|y|^2-k)^{\frac{1}{2}}\phi_r^{\frac{1}{2}}dx\right]^2\\\ns &\leq
4\varepsilon_2\displaystyle\int_\Omega|\nabla y|^2\phi_r
dx+C_1\varepsilon_2|y|^2_{H^1(\Omega;
\dbR^n)}+C_1\varepsilon_2^{-1}|y|^4_{L^2(\Omega; \dbR^n)}.
\end{array}
\end{eqnarray}
Therefore, substituting (\ref{20}) into (\ref{17}) and (\ref{45})
respectively, we see that
 $$
\begin{array}{ll}
\displaystyle\int_\Omega|f||y|\phi_r dx \leq\3n&
\left|f\right|_{L^{\frac{4\theta}{\theta+6}}(\Omega; \dbR^n)}\left[
4\varepsilon_2\displaystyle\int_\Omega|\nabla y|^2\phi_r
dx+C_1\varepsilon_2|y|^2_{H^1(\Omega;
\dbR^n)}+C_1\varepsilon_2^{-1}|y|^4_{L^2(\Omega; \dbR^n)}+1
\right]\\\ns
&+C_1\left|f\right|_{L^{\frac{\theta}{2}}(\Omega;
\dbR^n)}|y|^2_{H^1(\Omega; \dbR^n)}
\end{array}
 $$
and
\begin{eqnarray*}
&&\displaystyle \int_{\Omega}\left(\sum_{i,j=1}^n|C^{i
j}|^2+\sum_{i=1}^n|D^i|+1\right)|y|^2\phi_r dx\\\ns
&&\leq
4C_1L\varepsilon_2\displaystyle\int_\Omega|\nabla y|^2\phi_r
dx+C_1L(1+C_1\varepsilon_2)|y|^2_{H^1(\Omega;
\dbR^n)}+C_1^2L\varepsilon_2^{-1}|y|^4_{L^2(\Omega; \dbR^n)}.
\end{eqnarray*}

Combining the above inequalities with (\ref{14}) and taking
$\varepsilon_2$ sufficiently small such that
$$\left(4\left|f\right|_{L^{\frac{4\theta}{\theta+6}}(\Omega;
\dbR^n)}\varepsilon_2+4C_1L\varepsilon_2\right)C
<\displaystyle\frac{1}{2},$$ where $C$ and $C_1$ are the constants
appeared in (\ref{14}) and (\ref{45}) respectively, we arrive at
$$
\displaystyle\int_\Omega\left(|\nabla
y|^2\phi_r+|\nabla\phi_r|^2\right)dx \leq C_2,
$$
here and hereafter $C_2$ is a constant depending on $C$, $C_1$, $L$,
$|f|_{L^{\frac{\theta}{2}}(\Omega; \dbR^n)}$ and $|y|_{H^1(\Omega;
\dbR^n)}$, independent of $r$. Since $\phi_r\in H^1_0(\Omega)$, by
the definition of $\phi_r$, letting $r\rightarrow +\infty$ in the
above inequality, for any fixed $k\geq
\esssup\limits_{\Gamma}|y|^2$, we obtain that
\begin{equation}\label{27}
\displaystyle\int_\Omega|\nabla y|^2
(|y|^2-k)_+dx+\displaystyle\int_\Omega[(|y|^2-k)_+]^2
dx+\displaystyle\int_{|y|^2> k}|\nabla(|y|^2)|^2dx\leq C_2.
\end{equation}

Finally, we construct a sequence of inequalities with respect to
$A_k$, where $A_k=\{x\in \Omega; |y(x)|^2> k\}$. Again, by
(\ref{14}), we get that
\begin{eqnarray*}
\begin{array}{rl}
&\displaystyle\int_{\Omega}\left(|\nabla y|^2\phi_r+|\nabla
\phi_r|^2\right]dx\\\ns
&\leq
C\left[\displaystyle\int_{\Omega}|f||y|(|y|^2-k)_+dx+\displaystyle\int_\Omega\displaystyle\sum_{i=1}^n\Big(\displaystyle\sum_{j=1}^n
|C^{i j}|^2+|D^i|+1\Big)|y|^2(|y|^2-k)_+dx\right].
\end{array}
\end{eqnarray*}
Letting $r\rightarrow +\infty$ in the above inequality, for any
$\varepsilon_3>0$, by the H\"{o}lder inequality and Lemma \ref{31},
we see that
\begin{eqnarray}\label{24}
\begin{array}{rl}
&\displaystyle\int_{A_k}|\nabla y|^2
(|y|^2-k)dx+\displaystyle\int_{A_k}|\nabla |y|^2|^2 dx\\\ns &\leq
C\displaystyle\int_{A_k}\Big[|f|+\displaystyle\sum_{i=1}^n\Big(\displaystyle\sum_{j=1}^n
|C^{i j}|^2+|D^i|+1\Big)\Big](|y|^4+1)dx\\\ns &\leq
C(\left|f\right|_{L^{\frac{\theta}{2}}(\Omega; \dbR^n)}+L)
\left(\left||y|^2-k\right|_{L^{\frac{2\theta}{\theta-2}}(A_k)}^2
+k^2|A_k|^{1-\frac{2}{\theta}}\right)+C(\left|f\right|_{L^{\frac{\theta}{2}}
(\Omega; \dbR^n)}+L)|A_k|^{1-\frac{2}{\theta}}\\\ns &\leq
C\left(1+L+\left|f\right|_{L^{\frac{\theta}{2}}(\Omega;
\dbR^n)}\right)\left(\varepsilon_3|\nabla|y|^2|^2_{L^2(A_k)}+C(\varepsilon_3)||y|^2-k|_{L^2(A_k)}^2\right)\\\ns
&\quad+C\left(1+L+\left|f\right|_{L^{\frac{\theta}{2}}(\Omega;
\dbR^n)}\right)k^2|A_k|^{1-\frac{2}{\theta}}.
\end{array}
\end{eqnarray}
Denote $v=|y|^2$ and take $\varepsilon_3$ to be sufficiently small,
then by (\ref{27}) and (\ref{24}), one derives that
\begin{equation}\label{25}
\displaystyle\int_{A_k}|\nabla v|^2 dx\leq
C_3\displaystyle\int_{A_k}|v-k|^2dx+C_3k^2|A_k|^{1-\frac{2}{\theta}},
\end{equation}
where $C_3$ denotes a constant only depending only on $C$, $L$ and
$\left|f\right|_{L^{\frac{\theta}{2}}(\Omega; \dbR^n)}$.

By Lemma \ref{32}, we take
$$
m_0=2,\ l_0=2,\ \sigma=2,\
\varepsilon_0=\frac{2}{m}-\frac{2}{\theta},\ \gamma=C_3.
$$
Then it follows that
\begin{eqnarray}\label{26}
\begin{array}{rl}
&\esssup\limits_\Omega |y|\leq C\left(m,\ n,\ \theta,\ \Omega, \
\rho,\ |a^{i j}|_{L^\infty(\Omega)},\ L,\
\left|\displaystyle\frac{\det B^{i j}}{\det B}\right|_{W^{1, \infty}(\Omega)},\right.\\
&\left.\quad\quad\quad\quad\quad\quad\quad\quad
|f|_{L^\frac{\theta}{2}(\Omega; \dbR^n)},\ |y|_{L^2(\Omega;
\dbR^n)},\ \esssup\limits_{\Gamma}|y|\right).
\end{array}
\end{eqnarray}
Since $y$ is the weak solution, by Lemma \ref{2l1}, we have that
$$
|y|_{L^2(\Omega; \dbR^n)}\leq C(m,\ n,\ \theta,\ \Omega,\ \rho,\
|a^{i j}|_{{L^\infty}(\Omega)},\ L)(|f|_{L^\frac{\theta}{2}(\Omega;
\dbR^n)}+|g|_{H^1(\Omega; \dbR^n)}).
$$
This, combined with (\ref{26}), yields the desired conclusion in
Theorem \ref{5}.\endpf

\section{Proof of Theorem \ref{64}}\label{s5}

Now, let us prove our second main result, i.e.,
Theorem \ref{64}.\\

\noindent {\bf Proof of Theorem \ref{64}.} The main idea is the same
as that in the proof of Theorem \ref{5}. First, for any weak
solution $y=(y^1, y^2, \cdots, y^n)^\top$ to system $(\ref{61})$, We
choose $\varphi=(\varphi^1, \varphi^2, \cdots, \varphi^n)^\top\in
H^1_0(\Omega; \dbR^n)$ as a test function, where
$\varphi^i=\displaystyle\sum_{l=1}^{n}E^{i l}y^l\zeta_r$, and $E^{i
j}$ $(i, j=1,2, \cdots, n)$ are given by assumption {\bf (H)}, while
$\zeta_r$ is a suitable function to be specified later. By
Definition \ref{80}, it follows that
$$\begin{array}{rl} &\displaystyle\sum_{i,
j=1}^{n}\displaystyle\sum_{p,
q=1}^{m}\displaystyle\sum_{l=1}^{n}\displaystyle\int_{\Omega} a^{i
j}_{p q}y^j_{x_p}(E^{i l}y^l\zeta_r)_{x_q}dx+\displaystyle\sum_{i,
j, l=1}^{n}\displaystyle\int_{\Omega}C^{i j}\cdot\nabla y^jE^{i
l}y^l\zeta_rdx+\displaystyle\sum_{i,
l=1}^{n}\displaystyle\int_{\Omega}D^i\cdot yE^{i l}y^l\zeta_r
dx\\\ns
&=\displaystyle\sum_{i,
l=1}^{n}\displaystyle\int_{\Omega}f^iE^{i l}y^l\zeta_rdx.
\end{array}$$
This implies that
\begin{eqnarray*}
&&\displaystyle\sum_{i, j=1}^{n}\displaystyle\sum_{p,
q=1}^{m}\displaystyle\sum_{l=1}^{n}\displaystyle\int_{\Omega} \left[
a^{i j}_{p q}E^{i l}y^j_{x_p}y^l_{x_q}\zeta_r+a^{i j}_{p q}E^{i
l}y^j_{x_p}y^l(\zeta_r)_{x_q}+a^{i j}_{p q}(E^{i
l})_{x_q}y^j_{x_p}y^l\zeta_r \right]dx\\\ns
&&\q+\displaystyle\sum_{i, j, l=1}^{n}\displaystyle\int_{\Omega}C^{i
j}\cdot\nabla y^jE^{i l}y^l\zeta_rdx+\displaystyle\sum_{i,
l=1}^{n}\displaystyle\int_{\Omega}D^i\cdot yE^{i l}y^l\zeta_r
dx\\\ns
&&=\displaystyle\sum_{i,
l=1}^{n}\displaystyle\int_{\Omega}f^iE^{i l}y^l\zeta_rdx.
\end{eqnarray*}
Therefore,
\begin{eqnarray}\label{66}
\begin{array}{rl}
&\displaystyle\sum_{i, j=1}^{n}\displaystyle\sum_{p,
q=1}^{m}\displaystyle\sum_{l=1}^{n}\displaystyle\int_{\Omega} \left[
a^{i j}_{p q}E^{i l}y^j_{x_p}y^l_{x_q}\zeta_r+a^{i j}_{p q}E^{i
l}y^j_{x_p}y^l(\zeta_r)_{x_q}\right]dx\\\ns
&\leq C_4
\displaystyle\int_{\Omega}\left[|f||y|\zeta_r+\left(1+\sum_{i,j=1}^n|C^{i
j}|\right)|\nabla y||y|\zeta_r+\sum_{i=1}^n|D^i||y|^2\zeta_r
\right]dx,
\end{array}
\end{eqnarray}
here and hereafter $C_4$ denotes a constant depending only on $n$,
$m$, $\rho$, $|a^{i j}_{p q}|_{L^\infty(\Omega)}$ and $|E^{i
j}|_{W^{1, \infty}(\Omega)}$ $(i, j=1, \cdots, n; p, q=1, \cdots,
m)$.

Next, we estimate the terms in the left side of (\ref{66}). For this
purpose,  by (\ref{2e7}), condition 2) in assumption {\bf (H)} and
the Cramer rule, we see that for any $p, q=1, \cdots, m$, functions
$E^{i j}$ $(i, j=1, \cdots, n)$ (given by assumption {\bf (H)})
satisfy $ \displaystyle\sum_{l=1}^n a^{l i}_{p q}E^{l j}=f_{p q}h^{i
j}$. In particular, by $h^{1 1}=1$, we see that $f_{p
q}=\displaystyle\sum_{l=1}^{n}a^{l 1}_{p q}E^{l 1}.$ Therefore,
 \begin{equation}\label{81}
\displaystyle\sum_{l=1}^{n}a^{l i}_{p q}E^{l j}=h^{i
j}\displaystyle\sum_{l=1}^{n}a^{l 1}_{p q}E^{l 1}.
 \end{equation}
Hence,
\begin{eqnarray}\label{69}
\begin{array}{rl}
&\displaystyle\sum_{i, j=1}^{n}\displaystyle\sum_{p,
q=1}^{m}\displaystyle\sum_{l=1}^{n}\displaystyle\int_{\Omega}a^{i
j}_{p q}E^{i l}y^j_{x_p}y^l(\zeta_r)_{x_q}dx=\displaystyle\sum_{p,
q=1}^{m}\displaystyle\sum_{j,
l=1}^{n}\displaystyle\int_{\Omega}\left(\displaystyle\sum_{i=1}^{n}
a^{i 1}_{p q}E^{i 1}\right)h^{j l}y^j_{x_p}y^l(\zeta_r)_{x_q}dx\\\ns
&=\displaystyle\sum_{p,
q=1}^{m}\displaystyle\int_{\Omega}\left(\displaystyle\sum_{i=1}^{n}
a^{i 1}_{p q}E^{i
1}\right)\left[\displaystyle\frac{1}{2}\displaystyle\sum_{j=1}^{n}
h^{j j}(y^j)^2+\displaystyle\sum_{j, l\in\{1, 2, \cdots, n\},\
j<l}h^{j l}y^jy^l\right]_{x_p}(\zeta_r)_{x_q}dx\\\ns &\quad
-\displaystyle\sum_{p,
q=1}^{m}\displaystyle\int_{\Omega}\displaystyle\frac{1}{2}\left(\displaystyle\sum_{i=1}^{n}
a^{i 1}_{p q}E^{i 1}\right)\left[\displaystyle\sum_{j, l=1}^{n}(h^{j
l})_{x_p}y^jy^l\right](\zeta_r)_{x_q}dx.
\end{array}
\end{eqnarray}

On the other hand, by condition 4) in assumption {\bf (H)} and
noting (\ref{81}), it is easy to see that
$$M=\left(
\begin{array}{lcccccl}
\displaystyle\sum_{l=1}^{n}a_{1 1}^{l 1}E^{l 1}& \cdots &
\displaystyle\sum_{l=1}^{n}a_{1 m}^{l 1}E^{l 1} &
 \cdots &
\displaystyle\sum_{l=1}^{n}a_{1 1}^{l 1}E^{l n}& \cdots &
\displaystyle\sum_{l=1}^{n}a_{1 m}^{l 1}E^{l n}\\[2mm]
\vdots & \vdots & \vdots & \vdots  & \vdots & \vdots & \vdots \\
\displaystyle\sum_{l=1}^{n}a_{m 1}^{l 1}E^{l 1}& \cdots &
\displaystyle\sum_{l=1}^{n}a_{m m}^{l 1}E^{l 1} &
 \cdots &
\displaystyle\sum_{l=1}^{n}a_{m 1}^{l 1}E^{l n}& \cdots &
\displaystyle\sum_{l=1}^{n}a_{m m}^{l 1}E^{l n}\\[2mm]
\vdots & \vdots & \vdots & \vdots  & \vdots & \vdots & \vdots \\
\displaystyle\sum_{l=1}^{n}a_{1 1}^{l n}E^{l 1}& \cdots &
\displaystyle\sum_{l=1}^{n}a_{1 m}^{l n}E^{l 1} &
 \cdots &
\displaystyle\sum_{l=1}^{n}a_{1 1}^{l n}E^{l n}& \cdots &
\displaystyle\sum_{l=1}^{n}a_{1 m}^{l n}E^{l n}\\[2mm]
\vdots & \vdots & \vdots & \vdots & \vdots & \vdots & \vdots \\
\displaystyle\sum_{l=1}^{n}a_{m 1}^{l n}E^{l 1}& \cdots &
\displaystyle\sum_{l=1}^{n}a_{m m}^{l n}E^{l 1} &
 \cdots &
\displaystyle\sum_{l=1}^{n}a_{m 1}^{l n}E^{l n}& \cdots &
\displaystyle\sum_{l=1}^{n}a_{m m}^{l n}E^{l n}
\end{array}
\right)_{nm\times nm}.
$$
Therefore,
\begin{eqnarray}\label{67}
\begin{array}{rl}
&\displaystyle\sum_{i, j=1}^{n}\displaystyle\sum_{p,
q=1}^{m}\displaystyle\sum_{l=1}^{n}\displaystyle\int_{\Omega}  a^{i
j}_{p q}E^{i l}y^j_{x_p}y^l_{x_q}\zeta_r dx \geq \rho_3
\displaystyle\int_{\Omega} |\nabla y|^2\zeta_r dx.
\end{array}
\end{eqnarray}

Combining (\ref{69}) and (\ref{67}) with (\ref{66}), we have
\begin{eqnarray}\label{70}
\begin{array}{rl}
&\displaystyle\int_{\Omega} \left\{|\nabla
y|^2\zeta_r+\displaystyle\sum_{p,
q=1}^{m}\left(\displaystyle\sum_{i=1}^{n} a^{i 1}_{p q}E^{i
1}\right)\psi_{x_p}(\zeta_r)_{x_q}\right\}dx\\\ns &\leq C_5
\displaystyle\int_{\Omega}\left[|f||y|\zeta_r+\left(1+\sum_{i,j=1}^n|C^{i
j}|\right)|\nabla
y||y|\zeta_r+\sum_{i=1}^n|D^i||y|^2\zeta_r+|y|^2|\nabla\zeta_r|
\right]dx,
\end{array}
\end{eqnarray}
where $\psi=\displaystyle\sum_{j, l=1}^nh^{j l}y^jy^l$ and $C_5$
depends only on $C_4$, $\rho_3$ and $|h^{i j}|_{W^{1,
\infty}(\Omega)}$ $(i, j=1, 2, \cdots, n)$.

For $s,\ r>0$ and $k>\sup\limits_{\Gamma}\psi^s$, denote
$$A_k=\left\{x\in \Omega\ |\ \psi^s(x)>k \right\},\ \ \ \ A_k^r=\left\{x\in \Omega\ |\ k<\psi^s(x)<k+r \right\}.$$
Moreover, we choose $\zeta_r=\min\{r, (\psi^s-k)_+\}$. Then, by
(\ref{70}), and using condition 3) in assumption {\bf (H)},  we
conclude that
\begin{eqnarray}\label{71}
\begin{array}{rl}
&\displaystyle\int_{A_k} |\nabla y|^2\zeta_r
dx+\displaystyle\int_{A^r_k} \psi^{s-1}|\nabla\psi|^2 dx\\\ns &\leq
C_6\displaystyle\int_{A_k}
\left[|f||y|\zeta_r+\left(1+\sum_{i,j=1}^n|C^{i
j}|^2+\sum_{i=1}^n|D^i|\right)|y|^2\zeta_r+|y|^2|\nabla\zeta_r|
\right]dx,
\end{array}
\end{eqnarray}
where $C_6$ denotes a constant depending only on $s$, $C_5$ and
$\rho_2$. Moreover, using condition 1) in assumption {\bf (H)}, it
follows that
 $$\psi\geq \rho_1|y|^2.
 $$

Now, let us estimate the right side of (\ref{71}). First, by
H\"{o}lder's inequality, for any $\varepsilon_4>0$, we have that
\begin{equation}\label{72}
\displaystyle\int_{A_k}|f||y|\zeta_r dx\leq
C_7\displaystyle\int_{A_k}|y|\zeta_r dx\leq
\varepsilon_4\displaystyle\int_{A_k}\zeta_r^{\frac{s+1}{s}}dx+\varepsilon_4^{-1}C_7
\displaystyle\int_{A_k}|y|^{s+1}dx,
\end{equation}
here and hereafter $C_7$ denotes a constant depending only on $C_6$,
$\Omega$, $\rho_1$, $|f|_{L^\infty(\Omega; \dbR^n)}$,
$\displaystyle\sum_{i,j=1}^n|C^{i j}|_{L^\infty(\Omega; \dbR^m)}$
and $\displaystyle\sum_{i=1}^n|D^i|_{L^\infty(\Omega; \dbR^n)}$.
Next,
\begin{eqnarray}\label{73}
\begin{array}{rl}
&\displaystyle\int_{A_k}\left(1+\sum_{i,j=1}^n|C^{i
j}|^2+\sum_{i=1}^n|D^i|\right)|y|^2\zeta_r dx\leq
C_7\displaystyle\int_{A_k}|y|^2\zeta_r dx\\\ns &\leq
\varepsilon_4\displaystyle\int_{A_k}\zeta_r^{\frac{s+1}{s}}dx+\varepsilon_4^{-1}C_7\displaystyle\int_{A_k}|y|^{2s+2}dx.
\end{array}
\end{eqnarray}
Further,
\begin{eqnarray}\label{74}
\begin{array}{rl}
\displaystyle\int_{A_k}|y|^2|\nabla\zeta_r| dx\3n&\ds\leq
C_7\displaystyle\int_{A_k^r}|y|^2\psi^{s-1}|\nabla\psi| dx \leq
\varepsilon_4\displaystyle\int_{A_k^r}\psi^{s-1}|\nabla\psi|^2dx+\varepsilon_4^{-1}C_7\displaystyle\int_{A_k^r}|y|^4\psi^{s-1}dx\\\ns
&\leq
\varepsilon_4\displaystyle\int_{A_k^r}\psi^{s-1}|\nabla\psi|^2dx
dx+\varepsilon_4^{-1}C_7\displaystyle\int_{A_k^r}|y|^{2s+2}dx.
\end{array}
\end{eqnarray}

On the other hand, by Poinc\'{a}re's inequality,
\begin{equation}\label{75}
\displaystyle\int_{\Omega}\zeta_r^{\frac{s+1}{s}}dx\leq C_7
\displaystyle\int_{\Omega}\left|\nabla(\zeta_r^{\frac{s+1}{2s}})\right|^2dx\leq
C_7\displaystyle\int_{A^r_k}\zeta_r^{\frac{1-s}{s}}\psi^{2s-2}|\nabla\psi|^2dx\leq
C_7\displaystyle\int_{A^r_k}\psi^{s-1}|\nabla\psi|^2dx.
\end{equation}
Therefore, by (\ref{71})--(\ref{75}), taking $\varepsilon_4$
sufficiently small, we get that
$$
\displaystyle\int_{A_k}|\nabla y|^2\zeta_r
dx+\displaystyle\int_{A^r_k}|\nabla (\psi^{\frac{s+1}{2}})|^2dx\leq
C_7\displaystyle\int_{A_k}(|y|^{2s+2}+|y|^{s+1})dx.
$$
Letting $r\rightarrow +\infty$ in the above inequality, we have that
\begin{equation}\label{76}
\displaystyle\int_{A_k}|\nabla y|^2(\psi^s-k)
dx+\displaystyle\int_{A_k}|\nabla (\psi^{\frac{s+1}{2}})|^2dx\leq
C_7\displaystyle\int_{A_k}(|y|^{2s+2}+|y|^{s+1})dx.
\end{equation}
Notice that for any given constant $s\leq
\displaystyle\frac{2}{m-2}$ (if $m\leq 2$, $s$ can be any positive
number), the right side of (\ref{76}) is finite.

Denote $\tilde{k}=k^{\frac{s+1}{2s}}$ and
$A_{\tilde{k}}=\{x\in\Omega\ |\
\psi^{\frac{s+1}{2}}>\tilde{k}\}=A_k$. Then by (\ref{76}), it
follows that
\begin{eqnarray}\label{77}
\begin{array}{rl}
&\displaystyle\int_{A_{\tilde{k}}}|\nabla(\psi^{\frac{s+1}{2}})|^2dx\leq
C_7\displaystyle\int_{A_{\tilde{k}}}\psi^{s+1}dx+C_7\displaystyle\int_{A_{\tilde{k}}}\psi^{\frac{s+1}{2}}dx\\\ns
&\leq
C_7\displaystyle\int_{A_{\tilde{k}}}(\psi^{\frac{s+1}{2}}-\tilde{k})^2dx+C_7\tilde{k}^2|A_{\tilde{k}}|+
C_7\displaystyle\int_{A_{\tilde{k}}}(\psi^{\frac{s+1}{2}}-\tilde{k})dx+C_7\tilde{k}|A_{\tilde{k}}|\\\ns
&\leq
C_7\displaystyle\int_{A_{\tilde{k}}}(\psi^{\frac{s+1}{2}}-\tilde{k})^2dx+C_7\tilde{k}^2|A_{\tilde{k}}|.
\end{array}
\end{eqnarray}
By Lemma \ref{32}, we take
$$
u=\psi^{\frac{s+1}{2}},\ \ m_0=\sigma=l_0=2,\ \ k=\tilde{k},\ \
\gamma=C_7,\ \  \varepsilon_0=\displaystyle\frac{2}{m}.
$$
Then, using also Lemma \ref{2l1}, it follows that
\begin{eqnarray*}
&&\esssup\limits_\Omega |y|\leq C\left(m,\ n,\ \Omega, \ \rho,\
\rho_1,\ \rho_2,\ \rho_3,\ |a^{i j}_{p q}|_{L^\infty(\Omega)},\
|C^{i j}|_{L^\infty(\Omega; \dbR^m)},\ |D^i|_{L^\infty(\Omega;
\dbR^n)},\ \right.\\\ns &&\left.\quad\quad\quad\quad\quad\quad\quad
\left|E^{i j}\right|_{W^{1, \infty}(\Omega)},\ |h^{i j}|_{W^{1,
\infty}(\Omega)}, |g|_{H^1(\Omega; \dbR^n)},\ |f|_{L^\infty(\Omega;
\dbR^n)}, \ \esssup\limits_{\Gamma}|y|\right).
\end{eqnarray*}
This completes the proof of Theorem \ref{64}.
\endpf

\section{An example}\label{s6}

In this section we give an example, in which the coefficients $a^{i
j}_{p q}$ $(i, j=1, 2, \cdots, n; p, q=1, 2, \cdots, m)$ of system
(\ref{61}) satisfy all of the assumptions in Theorem \ref{64}.

For any given functions $b^{i j}\in W^{1, \infty}(\Omega)$ and $g_{p
q}\in L^\infty(\Omega)$ $(i, j=1, 2, \cdots, n; p, q=1, 2, \cdots,
m)$ such that $g_{p q}>0$  and the following matrix is uniformly
positive definite:
$$G:=\left(
\begin{array}{lccl}
g_{1 1} & g_{1 2} & \cdots & g_{1 m} \\[2mm]
g_{2 1} & g_{2 2} & \cdots & g_{2 m} \\[2mm]
\vdots & \vdots & \vdots & \vdots \\
g_{m 1} & g_{m 2} & \cdots & g_{m m}
\end{array}
\right),
$$
we take
$$a^{i j}_{p q}\ =\ b^{i j}g_{p q}, \quad\quad h^{i j}=\left\{\begin{array}{rll}&1 &\mbox{ if }i=j\\
&0 &\mbox{ if }i\neq j \end{array}\right.,\quad\quad f_{p q}\ =\
\displaystyle\frac{L_{p q}}{(g_{p q})^{n-1}}.$$ Then it is easy to
check the following assertions:

\begin{enumerate}

\item[i)] Condition 1) in assumption {\bf (H)} holds, since $V=I_{n\times n}$;

\item[ii)] By the definition of $a^{i j}_{p q}$ and $f_{p q}$,
and $b^{i j}\in W^{1, \infty}(\Omega)$ $(i, j=1, 2, \cdots, n; p,
q=1, 2, \cdots, m)$, if (\ref{2e7}) holds, condition 2) in
assumption {\bf (H)} is satisfied;

\item[iii)] If $b^{1 1}\neq0$, $b^{i 1}=0$, $i=2, 3, \cdots, n$,
and $E^{11}:=\displaystyle\frac{f_{p q}}{L_{p
q}}\displaystyle\sum_{l=1}^{n} h^{l 1}v^{l 1}_{p q}>0$ in $\Omega$,
then by (\ref{63}), conditions 3) and 4) in assumption {\bf (H)}
hold true. Notice that by the definition of $E^{11}$, $E^{11}>0$ if
and only if
\begin{eqnarray}\label{79}\mbox{ det }\left(
\begin{array}{lccl}
b^{2 2} & b^{2 3} & \cdots & b^{2 n} \\[2mm]
b^{3 2} & b^{3 3} & \cdots & b^{3 n} \\[2mm]
\vdots & \vdots & \vdots & \vdots \\
b^{n 2} & b^{n 3} & \cdots & b^{n n}
\end{array}
\right)>0.
\end{eqnarray}
Moreover, under condition (\ref{79}), the hypothesis (\ref{2e7})
also holds;

\item[iv)] Condition (\ref{63}) is equivalent to that the
following matrix is uniformly positive definite:
$$K:=\left(
\begin{array}{cccc}
b^{1 1}G & \displaystyle\frac{1}{2}b^{1 2}G & \cdots & \displaystyle\frac{1}{2}b^{1 n}G \\[6mm]
\displaystyle\frac{1}{2}b^{1 2}G &  b^{2 2}G & \cdots & \displaystyle\frac{1}{2}(b^{2 n}+b^{n 2})G \\[6mm]
\displaystyle\frac{1}{2}b^{1 3}G & \displaystyle\frac{1}{2}(b^{3 2}+b^{2 3})G & \cdots & \displaystyle\frac{1}{2}(b^{3 n}+b^{n 3})G \\[6mm]
\vdots & \vdots & \vdots & \vdots \\[6mm]
\displaystyle\frac{1}{2}b^{1 n}G & \displaystyle\frac{1}{2}(b^{2
n}+b^{n 2})G & \cdots & b^{n n}G
\end{array}
\right)_{nm\times nm}.
$$
Notice that if for some constant $\rho_*>0$,
\begin{equation}\label{78}
b^{i i}\geq \rho_* \quad\quad\quad\mbox{ and }\quad\quad\quad b^{i
j}\leq n\rho_* \quad\quad\quad(i, j=1, 2, \cdots, n; i\neq j),
\end{equation} then the matrix $K$ is uniformly positive definite.
\end{enumerate}

By the above assertions i)--iv), suppose that the coefficients $a^{i
j}_{p q}$ $(i, j=1, 2, \cdots, n$; $p, q=1, 2, \cdots, m)$ of system
(\ref{61}) satisfy that $$a^{i j}_{p q}\ =\ b^{i j}g_{p q},$$ where
$b^{i j}\in W^{1, \infty}(\Omega)$, $g_{p q}\in L^\infty(\Omega)$,
$b^{11}\neq0$, $b^{i 1}=0$ $(i=2, \cdots, n)$, $g_{p q}>0$, and $G$
is uniformly positive definite, and (\ref{79})--(\ref{78}) are
satisfied. Then, by Theorem \ref{64}, we conclude the boundedness of
weak solutions to the corresponding system (\ref{61}).

\end{document}